\documentclass[10pt]{article}    

\usepackage{amsmath,amssymb,graphicx,a4wide}

\numberwithin{equation}{section}

\begin{document}           

\title{Spectral Properties of Numerical Differentiation}

\author{Maxim Dvornikov    
  \\
  Department of Physics, 
  P.O.~Box~35, 40014, University of Jyv\"askyl\"a, 
  Finland\footnote{E-mail: dvmaxim@cc.jyu.fi} ;
  \\
  IZMIRAN, 142190, Troitsk, Moscow region, 
  Russia\footnote{E-mail: maxdvo@izmiran.ru}
  }         

\date{}

\maketitle                 

\thispagestyle{empty}

\begin{abstract}
We study the numerical differentiation formulae for functions
given in grids with arbitrary number of nodes. We investigate the
case of the infinite number of points in the formulae for the
calculation of the first and the second derivatives. The spectra
of the corresponding weight coefficients sequences are obtained.
We examine the first derivative calculation of a function given in
odd-number points and analyze the spectra of the weight
coefficients sequences in the cases of both finite and infinite
number of nodes. We derive the one-sided approximation for the
first derivative and examine its spectral properties.
\end{abstract}

\begin{center}
  Mathematics Subject Classification: Primary 65D25; Secondary 65T50
\end{center}


\section*{Introduction}

The problem of numerical differentiation is a long-standing issue.
There are plenty of published works devoted to the generation of
finite difference formulae in both one and multi dimensional
lattices (see, e.g., Ref.~\cite{AbrSti64}). However, many of those
methods require preliminary construction of an interpolating
polynomial, and hence are very awkward. Moreover, the majority of
the previous techniques are valid in the case of a function given
in the limited number of nodes.

The finite difference formulae for the calculation of any order
derivative in a one dimensional grid with arbitrary spacing were
discussed in Refs.~\cite{For88,For98}. However, only recursion
relations for the weight coefficients have been derived. The explicit formulas for 
the derivatives calculation were recently derived in Ref.~\cite{Li05} on the basis 
of the generalized Vandermonde determinant.

It should be noted that the low order derivatives (the first and
the second ones) as well as equidistant lattices are of the major
importance in many problems of applied mathematics and physics.
The first and the second numerical derivatives in the
equidistant one dimensional grid were studied in
Ref.~\cite{Dvo07JCAAM}. The finite difference formulae for the central
derivatives of a function given on a lattice with arbitrary number
of elements have been derived in that work. It is important that
these formulae have been obtained in the explicit form. This
method enabled one to examine the spectral properties of
weight coefficients sequences as well as to analyze the accuracy
of the numerical differentiation.

In the present paper we continue to study the numerical
differentiation formulae for functions given in grids with
arbitrary number of nodes. On the basis of the results of
Ref.~\cite{Dvo07JCAAM} in Sec.~\ref{FSDINP} we investigate the case of
the infinite number of points in the formulae for the calculation
of the first and the second derivatives. The spectra of the
corresponding weight coefficients sequences are also obtained.
Then, in Sec.~\ref{ONP} we examine the first derivative
calculation of a function given in odd-number points. We also
analyze the spectra of the weight coefficients sequences in the
cases of both finite and infinite number of nodes. In
Sec.~\ref{OSD} we derive the one-sided approximation for the first
derivative and examine its spectral properties. It is worth
noticing that the derivations of the finite difference formulae in
all cases are performed for the arbitrary number of points.
Finally, in Sec.~\ref{CONCL} we resume our results.

\section{Spectral properties of the first and the second
derivatives for infinite number of
points}\label{FSDINP}

Let us study the function $f(x)$ given in the equidistant points $x_m$, 
$f(x_m)=f_m$, where $m=0,\dots,\pm n$. It was found in Ref.~\cite{Dvo07JCAAM} that 
the first and the second derivatives are approximated as
\begin{equation}
  \label{f'}
  f^{\prime}(0)\approx 
  \frac{1}{2h}
  \sum_{m=1}^{n}\alpha_{m}^{(1)}(n)(f_{m}-f_{-m}).
\end{equation}
and
\begin{equation}
  \label{f''}
  f^{\prime\prime}(0)\approx 
  \frac{1}{h^{2}}
  \sum_{m=1}^{n}\alpha_{m}^{(2)}(n)(f_{m}-2f(0)+f_{-m}),
\end{equation}
where $h$ is the distance between nodes. The coefficients $\alpha_{m}^{(1)}(n)$ and 
$\alpha_{m}^{(2)}(n)$ can be calculated explicitly for arbitrary $n$ (see 
Ref.~\cite{Dvo07JCAAM}).

The spectral properties of the sequences $\alpha_{m}^{(1)}(n)$ and
$\alpha_{m}^{(2)}(n)$ in Eqs.~\eqref{f'} and~\eqref{f''} in the case of finite 
number of
interpolation points were carefully examined in
Ref.~\cite{Dvo07JCAAM}. We found out that the more points we involved
in the sequence $\alpha_{m}^{(1)}(n)$ the more close to linear the
corresponding spectrum was. Thus the considered sequence produces
more accurate first derivative of a function in the case of great
number of points. As for the sequence $\alpha_{m}^{(2)}(n)$, it
was also shown in Ref.~\cite{Dvo07JCAAM} that its spectrum approached
to parabola if $n>1$. We expect that the corresponding spectra
will be exactly linear and parabolic ones if $n\to\infty$.

Let us consider the spectral characteristics of the sequences
$\alpha_{m}^{(1)}(n)$ and $\alpha_{m}^{(2)}(n)$ in the case of
infinite number of interpolation points. First we remind the
result for the $\alpha_{m}^{(1)}(n)$ in the limit $n\to\infty$
(see Ref.~\cite{Dvo07JCAAM})
\begin{equation}
  \label{limalpha}
  \alpha_{m}^{(1)}=\lim_{n\to\infty}\alpha_{m}^{(1)}(n)=
  (-1)^{m+1}\frac{2}{m}.
\end{equation}
The Fourier transform of a function $f(x)$ can be presented in the
form (see, e.g., Ref.~\cite{MunWal00})
\begin{equation}\label{fourierinf}
  c(\omega)=h\sum_{x}e^{-i\omega x}f(x)=
  h\sum_{m=-\infty}^{+\infty}
  e^{-i\omega mh}f(mh).
\end{equation}
The inverse Fourier transformation is given by the following
expression:
\begin{equation*}
  f(x)=\int_{-\pi/h}^{\pi/h}
  \frac{d\omega}{2\pi}
  c(\omega)e^{i\omega x},
  \quad
  x=kh,
\end{equation*}
and has the cutoff at high frequencies, $\vert\omega\vert \leq \pi/h$.

\begin{figure}
  \centering
  \includegraphics[scale=.9]{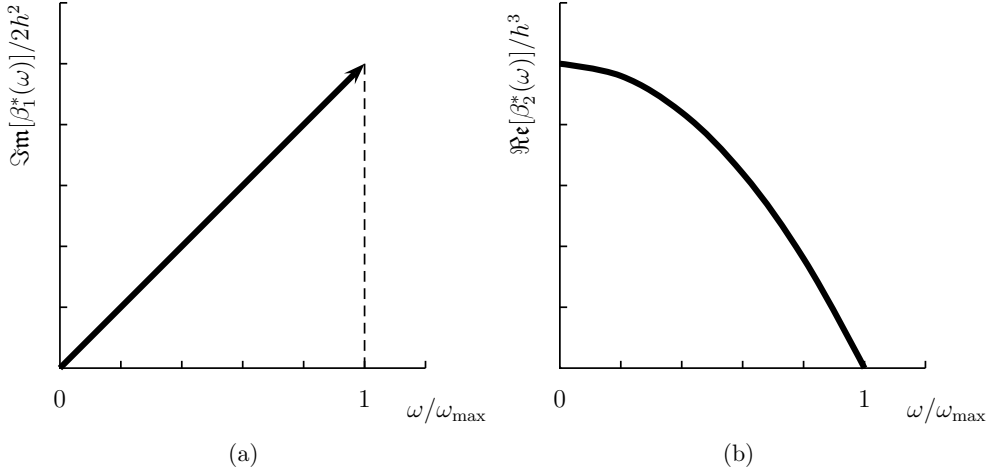}
  \caption{The spectra of differentiating filters, 
  (a) $\alpha_{m}^{(1)}$ and (b) $\alpha_{m}^{(2)}$, 
  in the case of infinite number of points.}
  \label{spbeta12}
\end{figure}

Now we can calculate the spectrum of the sequence $\alpha_{m}^{(1)}$,
\begin{equation}\label{beta1inf}
  \beta_1(\omega)= 
  h\sum_{\substack{m=-\infty\\ m\not=0}}^{+\infty}
  e^{-i\omega mh}\alpha_{m}^{(1)}
  = 
  -4ih\sum_{m=1}^{\infty}
  (-1)^{m-1}
  \frac{\sin(m\omega h)}{m}=
  -2i\omega h^2,
\end{equation}
where we use Eqs.~\eqref{limalpha} and \eqref{fourierinf}. Note
that Eq.~\eqref{beta1inf} is valid if $0\leq\omega<\pi/h$. The first derivative of 
the function can be expressed via the
spectra $\beta_1(\omega)$ and $c(\omega)$,
\begin{equation}\label{f'inf}
  f'(x)=
  \frac{1}{2h}
  \int_{-\pi/h}^{\pi/h}
  \frac{d\omega}{2\pi h}
  \beta^{*}_1(\omega)c(\omega)
  e^{i\omega x},
  \quad
  x=kh.
\end{equation}
Using the result for the calculation of $\beta_1(\omega)$
presented in Eq.~\eqref{beta1inf} we readily find that
\begin{equation}\label{f'inffin}
  f'(x)=
  \int_{-\pi/h}^{\pi/h}
  \frac{d\omega}{2\pi}
  (i\omega)c(\omega)
  e^{i\omega x}.
\end{equation}
Eq.~\eqref{f'inffin} shows that the first derivative calculation
with help of the sequence $\alpha_{m}^{(1)}$ gives the exact value
of the derivative in the case of infinite number of interpolation
points for all frequencies except
$\omega_{\mathrm{max}}=\pi/h$. The fact that the first derivative
computation does not give correct results at
$\omega=\omega_{\mathrm{max}}$ also follows from
Fig.~\ref{spbeta12}(a). However, it can be verified directly with help
of Eq.~\eqref{limalpha} for the function
$f_m=(-1)^m=\cos(\omega_{\mathrm{max}}mh)$.

In Ref.~\cite{Dvo07JCAAM} we showed that the use of the
coefficients $\alpha_{m}^{(1)}$ give exact value for the first
derivative of the function $y(x)=\sin(\omega_\mathrm{max}x/2)$.
However, Eq.~\eqref{f'inffin} [see also Fig.~\ref{spbeta12}(a)]
indicates that this method will give the correct results not only
for $\omega=\omega_\mathrm{max}/2$ but also for all frequencies
$\omega<\omega_\mathrm{max}$.

The second derivative calculation in the case of infinite number
of interpolation points can be analyzed in the similar manner as
we have done it for the first derivative. The explicit form of the
sequence $\alpha_{m}^{(2)}$ is
\begin{equation*}
  \alpha_{m}^{(2)}=\lim_{n\to\infty}\alpha_{m}^{(2)}(n)=
  (-1)^{m+1}\frac{2}{m^2}.
\end{equation*}
For the spectrum of the sequence $\alpha_{m}^{(2)}$ we obtain
\begin{equation}
  \label{beta2inf}
  \beta_2(\omega)=
  -\omega^2 h^3+
  \frac{\pi^2}{3}h.
\end{equation}
It should be noted that Eq.~\eqref{beta2inf} is valid for all
frequencies $0\leq\omega\leq\pi/h$. The expression for the second
derivative takes the form
\begin{equation}
  \label{f''inf}
  f''(x)= 
  \frac{1}{h^2}
  \int_{-\pi/h}^{\pi/h}
  \frac{d\omega}{2\pi h}
  \left[
  \beta^{*}_2(\omega)-\beta^{*}_2(0)
  \right]
  c(\omega)
  e^{i\omega x}
  = 
  \int_{-\pi/h}^{\pi/h}
  \frac{d\omega}{2\pi}
  (-\omega^2)c(\omega)
  e^{i\omega x},
  \quad
  x=kh.
\end{equation}
Eq.~\eqref{f''inf} demonstrates that the computation of the second
derivative with the use of the sequence $\alpha_{m}^{(2)}$ gives
the exact results in the case of infinite number of interpolation
points for all frequencies even including the maximal one. The
spectrum $\beta_2(\omega)$ is depicted in Fig.~\ref{spbeta12}(b).
%

\section{The first derivative computation of a
function given in odd-number points}\label{ONP}

In this section we discuss the calculation of the first derivative
in the case of a function given in odd-number nodes. Then we
discuss the spectral properties of the derived weight coefficients
sequences in the case of both finite and infinite number of nodes.

It follows from Fig.~\ref{spbeta12}(a) that the computation of the
first derivative gives unsatisfactory results at high frequencies
near $\omega_\mathrm{max}$. In order to introduce the numerical
differentiation of such rapidly oscillating functions we
consider the modified sequence
\begin{equation}\label{alpha1/2}
  \alpha_{2m+1}^{(1/2)}(n)=
  \frac{1}{(2m+1)\pi_m^{(1/2)}(n)},
  \quad
  m=0,\dots,n-1,
\end{equation}
where
\begin{equation*}
  \pi_m^{(1/2)}(n)=\prod_{
  \substack{k=0
  \\
  k{\not=}m}}^{n-1}
  \left(
  1-\frac{(2m+1)^{2}}{(2k+1)^{2}}
  \right),
\end{equation*}
and $\alpha_{2m}^{(1/2)}(n)=0$.

The coefficients in Eq.~\eqref{alpha1/2} can be formally derived
if we consider the first derivative calculation of a function
given in the odd-number points only
\begin{equation}\label{f'1/2}
  f'(0)\approx
  \frac{1}{2h}
  \sum_{m=0}^{n-1}\alpha_{2m+1}^{(1/2)}(n)(f_{2m+1}-f_{-2m-1}).
\end{equation}
Note that originally the function $f(x)$ was given in $2n+1$
points.

It is worth noticing that the results for the computation of the
weights with help of Eq.~\eqref{alpha1/2} in some particular cases
(namely for $n=3,5,7$ and $9$) coincide with those presented in
Ref.~\cite{For88} for the centered approximations at a 'half-way'
point. However, the method for the central derivatives calculation
elaborated in our paper enables one to get the expressions for the
weight coefficients in the explicit form for any number of nodes.

\begin{figure}
  \centering
  \includegraphics{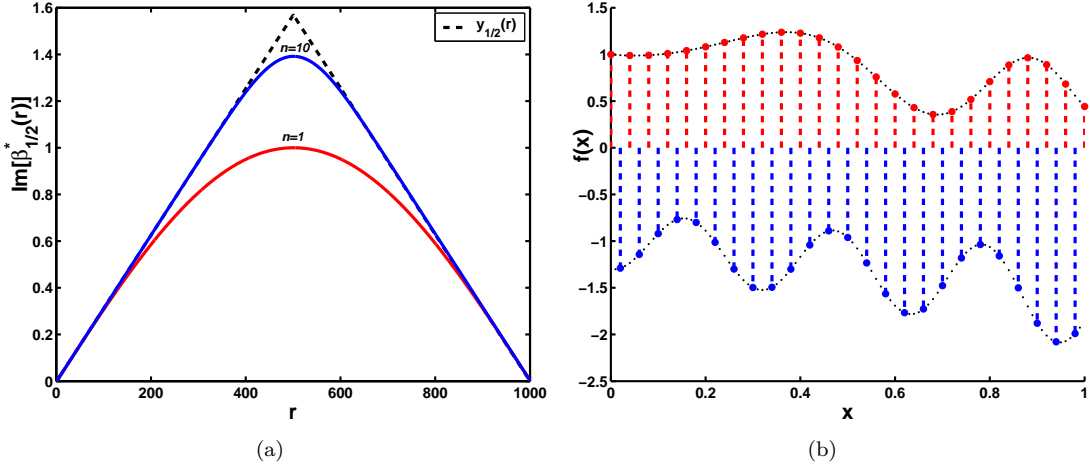}
  \caption{(a) The spectra of various sequences $\alpha_{m}^{(1/2)}(n)$
  at $N=2000$. (b) The example of the use of the sequence $\alpha^{(1/2)}_m(n)$.
  A rapidly oscillating function corresponds to the dashed lines,
  its envelope functions are shown by the dotted lines.}
  \label{spbeta1/2}
\end{figure}

We consider the spectral properties of the obtained sequence
$\alpha_{m}^{(1/2)}(n)$. Using the technique developed in
Ref.~\cite{Dvo07JCAAM} one can compute the spectrum of the sequence in
question,
\[
  \beta_{1/2}(r)=\sum_{m=0}^{N-1}
  \alpha_{m}^{(1/2)}(n)\exp
  \left(
  -i\frac{2\pi}{N}mr
  \right).
\]
The spectra of the sequences $\alpha_{m}^{(1/2)}(n)$ are depicted
in Fig.~\ref{spbeta1/2}(a) for the various values of $n$ at $N=2000$.
The function $y_{1/2}(r)$ has the form
\[
  y_{1/2}(r)=
  \begin{cases}
    2\pi r/N, & \text{if $0\leq r\leq N/4$}, \\
    \pi-2\pi r/N, & \text{if $N/4\leq r\leq N/2$}.
  \end{cases}
\]
It follows from this figure that for $n=1$ the imaginary part of
the spectrum is the function $\sin(2\pi r/N)$. The linearity
condition is satisfied only in the vicinity of zero and $N/2$.
However at $n=10$ linearity condition remains valid for $r\lesssim
350$ and $r\gtrsim 650$.

We examine the details of the differentiation procedure performed
with help of the sequence $\alpha_{m}^{(1/2)}(n)$. If a
function is slowly varying, then Eq.~\eqref{f'1/2} gives the
approximate value of the first derivative. It also results from
Fig.~\ref{spbeta1/2}(a). However, if a function is rapidly
oscillating (e.g., at $\omega\lesssim\omega_\mathrm{max}$), we can
consider its upper and lower envelope functions. Therefore the use
of the sequence $\alpha_{m}^{(1/2)}(n)$ will produce the first
derivative of envelope functions since in Eq.~\eqref{f'1/2} we
calculate the sum in odd-number points only. Envelope functions
can be treated as "smooth" in this case. Fig.~\ref{spbeta1/2}(b)
schematically illustrates this process.

\begin{figure}
  \centering
  \includegraphics[scale=.8]{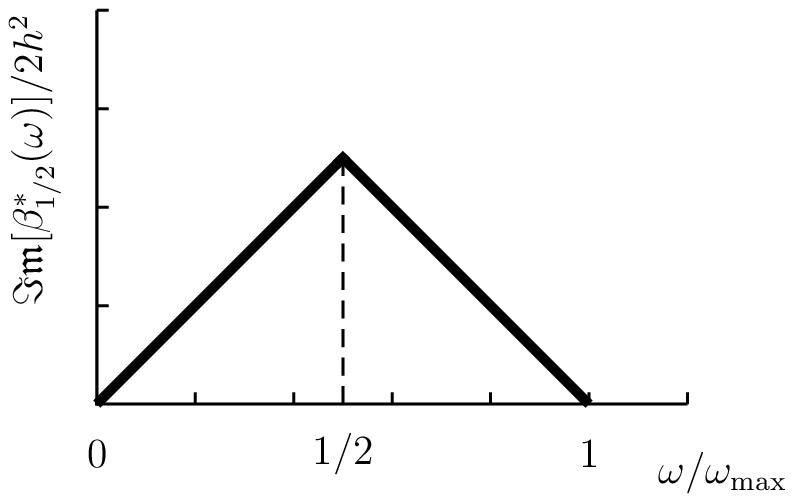}
  \caption{The spectrum of the sequence $\alpha_{m}^{(1/2)}$
  in the case of infinite number of points.}
  \label{spbeta1/2inf}
\end{figure}

Now let us discuss the case of infinite number of the
interpolation points. We can treat the sequence
$\alpha_{m}^{(1/2)}(n)$ in the similar way as it was done in
Ref.~\cite{Dvo07JCAAM}. Indeed, proceeding to the limit $n\to\infty$ in
Eq.\eqref{alpha1/2} we find that
\begin{equation}\label{alpha1/2lim}
  \alpha_{2m+1}^{(1/2)}=
  \lim_{n\to\infty}\alpha_{2m+1}^{(1/2)}(n)=
  (-1)^m
  \frac{4}{\pi(2m+1)^2}.
\end{equation}
In Eq.~\eqref{alpha1/2lim} we used the known value of infinite
product,
\[
  \prod_{k=0}^{\infty}
  \left(
  1-\frac{x^{2}}{(2k+1)^{2}}
  \right)
  =\cos\left(\frac{\pi x}{2}\right).
\]
With help of Eq.~\eqref{alpha1/2lim} it is possible to obtain the
spectrum of the sequence $\alpha_{m}^{(1/2)}$
\begin{equation}\label{beta1/2inf}
  \beta_{1/2}(\omega)=-2ih\times
  \begin{cases}
    \omega h, & \text{if $0\leq\omega\leq\pi/2h$}, \\
    (\pi-\omega h), & \text{if $\pi/2h\leq\omega\leq\pi/h$}.
  \end{cases}
\end{equation}
The imaginary part of the spectrum $\beta_{1/2}(\omega)$ is
presented in Fig.~\ref{spbeta1/2inf}.
As it follows from Eq.~\eqref{beta1/2inf} (see also
Fig.~\ref{spbeta1/2inf}) the sequence $\alpha_{m}^{(1/2)}$
performs the differentiation of a function in question if
$\omega\leq\omega_\mathrm{max}/2$, or its envelope functions if
$\omega\geq\omega_\mathrm{max}/2$. The of case
$\omega=\omega_\mathrm{max}/2$ was also considered in details in
Ref.~\cite{Dvo07JCAAM}.

\section{One-sided approximation for the first derivative}\label{OSD}

In this section we derive the weight coefficients for the
one-sided approximation of the first derivative and then we
analyze the spectral characteristics of the weight coefficients
sequence. It should be noted that the derivation of the weight
coefficients is analogous to case the of the central derivatives
which was carefully examined in Ref.~\cite{Dvo07JCAAM}.

Without restriction of generality we suppose that we
approximate the first derivative in the zero point. Let us
consider the function $f(x)$ given in the equidistant nodes
$x_{m}=mh>0$, where $m=0,\dots,n$, and $h$ is the constant value.
We can pass the interpolating polynomial of the $n$th power
through these points,
\begin{equation*}
  P_{n}(x)=\sum_{k=0}^{n}c_{k}x^{k}.
\end{equation*}
The values of the function in the nodes $x_{m}=mh$,
$f_{m}=f(x_{m})$, should coincide with the values of the
interpolating polynomial in these points,
\begin{equation}
  \label{fm}
  f_{m}=\sum_{k=0}^{n}c_{k}h^{k}m^{k}.
\end{equation}
In order to find the coefficients $c_{k}$, $k=0,\dots,n$, we
receive the system of inhomogeneous linear equations with the
given free terms $f_{m}$. It will be shown below that this system
has the single solution.

We will seek the solution of the system \eqref{fm} in the
following way:
\begin{equation*}
  c_{k}=\frac{1}{h^k}
  \sum_{m=0}^{n}f_{m}a_{m}^{(k)}(n),
\end{equation*}
where $a_{m}^{(k)}(n)$ are the undetermined coefficients
satisfying the condition,
\begin{equation}
  \label{a}
  \sum_{m=0}^{n}a_{m}^{(l)}(n)m^{k}=\delta_{lk},
  \quad l,k=0,\dots n.
\end{equation}
It is worth to be noted that, if we set $k=0$ and $l{\not=}0$ in
Eq.~\eqref{a}, we obtain the constraint which should be imposed on
the coefficients $a_{m}^{(l)}(n)$
\begin{equation}
  \label{acond}
  \sum_{m=0}^{n}a_{m}^{(l)}(n)=0,
  \quad l=1,\dots n.
\end{equation}
Analogous relation between the weight coefficients was established
in Ref.~\cite{For88}. In deriving of Eq.~\eqref{acond} (as well as
in all subsequent similar formulae) we suppose that $m^0=1$ if
$m=0$.

Let us resolve the system of equations \eqref{a} according to the
Cramer's rule
\begin{equation}
  \label{kram}
  a_{m}^{(l)}(n)=
  \frac{\Delta_{m}^{(l)}(n)}{\Delta_{0}(n)},
\end{equation}
where
\begin{equation}
  \label{delta0}
  \Delta_{0}(n)=
  \begin{vmatrix}
  1 & 1 & 1 & \dots & 1
  \\
  0 & 1 & 2 & \dots & n
  \\
  0 & 1 & 2^{2} & \dots & n^{2}
  \\
  \hdotsfor{5}
  \\
  0 & 1 & 2^n & \dots & n^n
  \\
  \end{vmatrix}=
  n!\prod_{1\leq i<j\leq n}(j-i){\not=}0,
\end{equation}
and
\begin{equation}
  \label{deltam}
  \Delta_{m}^{(l)}(n)=
  \begin{vmatrix}
  1 & 1 & 1 & \dots & 1 & 0 & 1 & \dots & 1
  \\
  0 & 1 & 2 & \dots & m-1 & 0 & m+1 & \dots & n
  \\
  \hdotsfor{9}
  \\
  0 & 1 & 2^{l} & \dots & (m-1)^{l} & 1
  & (m+1)^{l} & \dots & n^{l}
  \\
  \hdotsfor{9}
  \\
  0 & 1 & 2^{n} & \dots & (m-1)^{n} & 0
  & (m+1)^{n} & \dots & n^{n}
  \end{vmatrix}.
\end{equation}
In Eq.~\eqref{delta0} we use the formula for the calculation of
the Vandermonde determinant. From Eq.~\eqref{delta0} it follows
that the determinant of the system of equations \eqref{a} is not
equal to zero, i.e. the system of equations \eqref{fm} has the
single solution.

The most simple expression for $\Delta_{m}^{(l)}(n)$ is obtained
in the case of $l=1$ that corresponds to the calculation of the
first-order derivative
\begin{equation}
  \label{deltam1}
  \Delta_{m}^{(1)}(n)=(-1)^{m+1}
  \left(
  \frac{n!}{m}
  \right)^{2}
  \prod_{
  \substack{1\leq i<j\leq n
  \\
  i,j{\not=}m}}(j-i),
  \quad
  m=1,\dots,n.
\end{equation}
From Eq.~\eqref{kram} as well as taking into account
Eqs.~\eqref{delta0} and \eqref{deltam1} we get the expression for
the coefficients $a_{m}^{(1)}(n)$
\begin{equation}
  \label{am}
  a_{m}^{(1)}(n)=(-1)^{m+1}\frac{1}{m}\binom{n}{m},
  \quad
  m=1,\dots,n,
\end{equation}
where
\[
  \binom{n}{m}=\frac{n!}{m!(n-m)!},
\]
are the binomial coefficients. It is remarkable to note that the
coefficient $a_{1}^{(1)}(n)=n$. To simplify numerical calculations
(especially the analysis of the spectra of the derived sequences)
Eq.~\eqref{am} should be rewritten in the form
\begin{equation*}
  a_{m}^{(1)}(n)=\frac{1}{m p_m(n)},
  \quad
  m=1,\dots,n,
\end{equation*}
where
\[
  p_{m}(n)=\prod_{
  \substack{k=1
  \\
  k{\not=}m}}^{n}
  \left(
  1-\frac{m}{k}
  \right).
\]

In order to find the coefficient $a_{0}^{(1)}(n)$ we use
Eq.~\eqref{acond} rather than compute the determinant
\eqref{deltam}. Thus we receive the following expression for this
coefficient,
\begin{equation}\label{a0}
  a_{0}^{(1)}(n)=-\sum_{m=1}^{n}a_{m}^{(l)}(n)=
  -\sum_{m=1}^{n}(-1)^{m+1}\frac{1}{m}\binom{n}{m}=
  -\sum_{m=1}^{n}\frac{1}{m}.
\end{equation}
Eqs.~\eqref{am} and \eqref{a0} provide the weight coefficients for
the one-sided approximation of the first derivative of the
function $f(x)$ given in $n+1$ equidistant nodes,
\[
  f'(0)\approx
  \frac{1}{h}
  \sum_{m=0}^{n}f_{m}a_m^{(1)}(n).
\]
The results for the computation of the weights in some particular
cases (namely for $n=1,2,\dots,8$) coincide with those presented
in Ref.~\cite{For88}. However, the technique for derivatives
calculation developed in the present work allows one to obtain the
expressions for the weight coefficients in the explicit form for
any $n$.

Without derivation we mention that on the basis of
Eqs.~\eqref{kram}-\eqref{deltam} one can find the coefficients
$a_{m}^{(n)}(n)$ that correspond to the computation of the
$n$th-order derivative
\begin{equation}\label{an}
  a_{m}^{(n)}(n)=(-1)^{m+n}
  \frac{1}{n!}\binom{n}{m}.
  \quad
  m=0,\dots,n.
\end{equation}
It should be noted that Eq.~\eqref{an} is consistent with
Eq.~\eqref{acond}.

\begin{figure}
  \centering
  \includegraphics{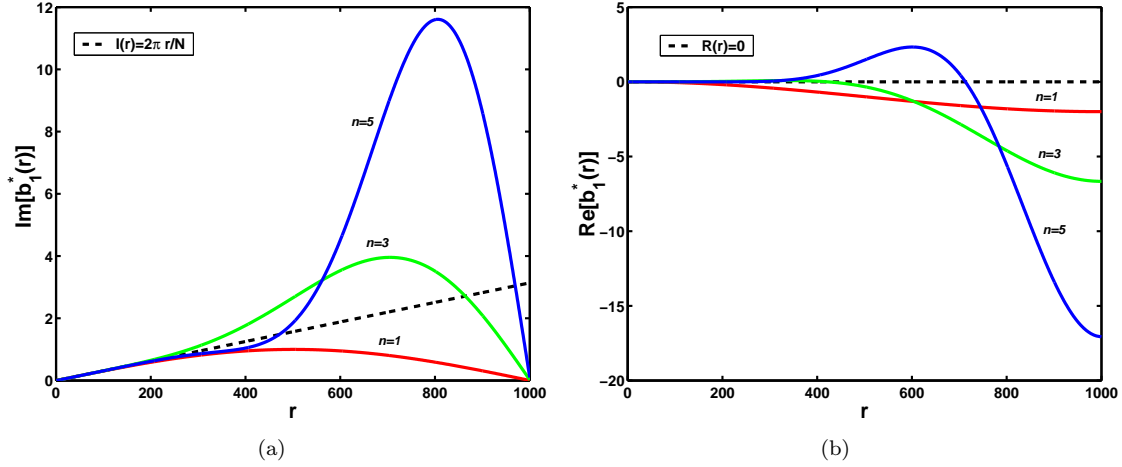}
  \caption{The imaginary parts (a) and the real parts (b) 
  of the spectra of various sequences $a_{m}^{(1)}(n)$ at $N=2000$.}
  \label{sponeside}
\end{figure}

Now we consider the spectral properties of the derived sequence
$a_{m}^{(1)}(n)$. Using the results of the previous section (see
also Ref.~\cite{Dvo07JCAAM}) we readily find the expression for the
spectrum of the considered sequence,
\[
  b_{1}(r)=\sum_{m=0}^{N-1}
  a_{m}^{(1)}(n)\exp
  \left(
  -i\frac{2\pi}{N}mr
  \right).
\]
The imaginary parts, $\Im\mathfrak{m}\big[b^*_1(r)\big]$, of the
spectra of the sequences $a_{m}^{(1)}(n)$ for the various values
of $n$ at $N=2000$ as well as the linearly growing sequence
$I(r)=2\pi r/N$ are presented in Fig.~\ref{sponeside}(a).
It follows form the this figure that the imaginary parts are close
to the linear sequence only in the vicinity of zero ($r\lesssim
200$) even at $n=5$. It points out that the one-sided
approximation of the first derivative has worse accuracy in
comparison with the central derivatives. This fact was also
mentioned in Ref.~\cite{DemMar63eng}. Therefore, the application of the
sequence $a_{m}^{(1)}(n)$ for the calculation of the one-sided
first derivative will give reliable results only for slowly
varying functions.

The spectrum $b_1(r)$ has not only imaginary part, but also
nonzero real part since $a_{0}^{(1)}(n){\not=}0$. The real parts,
$\Re\mathfrak{e}\big[b^*_1(r)\big]$, of the spectra of the
sequences $\alpha_{m}^{(1)}(n)$ for the various values of $n$ at
$N=2000$ as well as the constant sequence $R(r)=0$ are shown in
Fig.~\ref{sponeside}(b). 
It can be seen from this figure that the real parts of the spectra
are close to zero if $r\lesssim 100$ at $n=1$, and if $r\lesssim
300$ at $n=3$ and $n=5$. The deviation from zero is especially
great if $r\gtrsim 700$ at $n=3$, and if $r\gtrsim 500$ at $n=5$.
Such a behavior of the real parts of the spectra also reveals the
limited level of accuracy of the one-sided first derivative.

\section{Conclusion}\label{CONCL}

In conclusion we note that in our paper we have studied the
numerical differentiation formulae for functions given on grids
with arbitrary number of nodes. In Sec.~\ref{FSDINP} we have
investigated the case of the infinite number of points in the
formulae for the calculation of the first and the second
derivatives. The spectra of the corresponding weight coefficients
sequences have been obtained. It has been revealed that the
calculation of the first derivative with help of the derived
formulae gave reliable results for all spacial frequencies except
$\omega_\mathrm{max}$. As for the calculation of the second
derivative we have shown that the corresponding formulae were
valid for all spacial frequencies including
$\omega_\mathrm{max}$. In Sec.~\ref{ONP} we have examined the
first derivative calculation of a function given in odd-number
points. We have also analyzed the spectra of the weight
coefficients sequences in the cases of both finite and infinite
number of nodes. It has been found out that the obtained formulae
perform the differentiation of the considered function if
$\omega\leq\omega_\mathrm{max}/2$, and its envelope functions if
$\omega\geq\omega_\mathrm{max}/2$. In Sec.~\ref{OSD} we have
derived the one-sided approximation for the first derivative and
examined its spectral properties. The accuracy of the
one-sided first derivative has been discussed. On the basis of the
spectral properties of the weight coefficients sequences it has
been shown that the accuracy of the one-sided approximation for
the first derivative was essentially lower compared to the
computation of the central derivatives. This our result is in
agreement with previous works (see, e.g. Ref.~\cite{DemMar63eng}).
Nevertheless, the obtained one-sided first derivative formulae
could be of use in solving differential equations by means of
numerical methods. It is also possible to apply the elaborated technique of the 
numerical differentiation in construction of quantum field theory models of unified 
interactions.

\subsection*{Acknowledgments}

This research has been supported by the Academy of Finland under the 
contract No.~108875. The author is indebted to the Russian Science Support 
Foundation for a grant as well as to Sergey Dvornikov for
helpful discussions.


\begin{thebibliography}{7}

\bibitem{AbrSti64}
  M. Abramowitz and I.~A. Stegun,
  \textit{Handbook of Mathematical Functions},
  National Bureau of Standards, Washington D. C., 1964.
  
\bibitem{For88}
  B. Fornberg,
  Generation of Finite Difference Formulas on Arbitrary Spaced
  Grids,
  \textit{Math. Comp.}, 51(184), 699-706, (1988).  
  
\bibitem{For98}
  B. Fornberg,
  Calculation of Weights in Finite Difference Formulas,
  \textit{SIAM Rev.}, 40(3), 685-691, (1998).
  
\bibitem{Li05}
  J. Li,
  General Explicit Difference Formulas for Numerical Differentiation,
  \textit{J. Comp. \& Appl. Math.}, 183, 29-52, (2005).  
  
\bibitem{Dvo07JCAAM}
  M. Dvornikov,
  Formulae of Numerical Differentiation,
  \textit{JCAAM}, 5, 77-88, (2007) [e-print arXiv: math.NA/0306092].
  
\bibitem{MunWal00}
  G. M\"{u}nster and M. Walzl,
  Lattice Gauge Theory -- A Short Primer,
  in: \textit{Proceedings of the Summer School on Phenomenology of 
  Gauge Interactions} (D. Graudenz and V. Markushin, eds.),
  PSI, Villigen, Switzerland, 2000, pp.~127-160 [e-print arXiv: hep-lat/0012005]. 
  
\bibitem{DemMar63eng}
  B.~P. Demidovich and I.~A. Maron,
  \textit{Foundations of Computational Mathematics} (2nd ed.),
  Fiz. Mat. Lit., Moscow, 1963.

\end{thebibliography}
\end{document}